\documentclass[12pt,a4paper]{amsart}
\usepackage[utf8]{inputenc}	
\usepackage[T1]{fontenc}
\usepackage[in,headings]{fullpage}
\usepackage{amsbsy, amscd, amsfonts, amsmath, amsrefs, amssymb, amsthm}
\usepackage{graphicx}
\usepackage{float}
\usepackage{color}   
\usepackage[colorlinks=true,linkcolor=blue,citecolor=blue,urlcolor=blue]{hyperref}
\usepackage{comment}
\excludecomment{A}

\newtheorem{theorem}{Theorem}

\newtheorem{corollary}[theorem]{Corollary}
\newtheorem{proposition}[theorem]{Proposition}
\newtheorem{lemma}[theorem]{Lemma}

\theoremstyle{definition}
\newtheorem{definition}[theorem]{Definition}
\newtheorem{remark}[theorem]{Remark}

\theoremstyle{remark}

\newcommand{\C}{\mathbf{C}}
\newcommand{\Z}{\mathbf{Z}}

\newcommand{\R}{\mathbf{R}}

\renewcommand{\Re}{\mathop{\mathrm{Re}}\nolimits}
\renewcommand{\Im}{\mathop{\mathrm{Im}}\nolimits}
\newcommand{\Rzeta}{\mathop{\mathcal R }\nolimits}

\newfont{\cmbsy}{cmbsy10}
\newfont{\cmmib}{cmmib10}
\newcommand{\Orden}{\mathop{\hbox{\cmbsy O}}\nolimits}
\newcommand{\orden}{\mathop{\hbox{\cmmib o}}\nolimits}

\begin{document}

\title{Infinite Product of the Riemann auxiliary function.}
\author[Arias de Reyna]{J. Arias de Reyna}
\address{%
Universidad de Sevilla \\ 
Facultad de Matem\'aticas \\ 
c/Tarfia, sn \\ 
41012-Sevilla \\ 
Spain.} 




\email{arias@us.es, ariasdereyna1947@gmail.com}


\begin{abstract}
We obtain the product for the auxiliary function $\Rzeta(s)$ and study some related functions  as its phase $\omega(t)$ at the critical line. The function $\omega(t)$ determines the zeros of $\zeta(s)$ on the critical line. We study the influence of the zeros of $\Rzeta(s)$ on $\omega(t)$. Thus, the relationship between the zeros of $\Rzeta(s)$ and those of $\zeta(s)$ is determined.  
\end{abstract}

\maketitle
\section{Introduction}
The function  
\begin{displaymath}
\Rzeta(s)=\int_{0\swarrow1}\frac{x^{-s} e^{\pi i x^2}}{e^{\pi i x}-
e^{-\pi i x}}\,dx.
\end{displaymath}
was studied by Riemann and introduced by Siegel in \cite{Siegel}. It is an entire function. We show that it is of order 1, and establish its infinite product expansion. 

We apply this product expansion to study the phase $\omega(t)$ of $\Rzeta(\frac12+it)$, this is a real, real analytic function defined on $\R$ and such that $\Rzeta(\frac12+it)=e^{-i\omega(t)}g(t)$ where $g$ is also a real and real analytic function. Each point $t\in\R$ at which  $\cos(\vartheta(t)-\omega(t))=0$ determine a zero $\frac12+it$ of $\zeta(s)$ on the critical line. We show that 
$\omega(t)\approx2\pi N_r(t)$  where $N_r(t)$ is the number of zeros $\varrho=\beta+i\gamma$ of $\Rzeta(s)$  with $0<\gamma\le t$ and $\beta>1/2$. Therefore, if $\omega(t)=\orden(t)$, all except an infinitesimal proportion of zeros of $\zeta(s)$ will be on the critical line. 

Our computation of zeros (see \cite{A172}) of $\Rzeta(s)$  suggests that $N_r(t)\approx \frac{\vartheta(t)}{6\pi}$. If this were true $\vartheta(t)-\omega(t)\approx 
\frac{2}{3}\vartheta(t)$, and at least $\frac23$ of the zeros of $\zeta(s)$ will be on the critical line. The method of Levinson \cite{Lev}, has been refined by several authors  \cite{Con}, \cite{F}, \cite{P}, to obtain $\frac{5}{12}$ zeros of $\zeta(s)$ on the critical line. Therefore, the study of the zeros of $\Rzeta(s)$  has the potential  to improve these results. 

We give some simple Theorems about $\omega(t)$ and certain related functions.

\section{Some facts and notation on Theta functions.}
To determine the order of $\Rzeta(s)$ we will use the integral representation
given in \cite{A166}
\begin{theorem}
For $s\ne 0$
\begin{equation}\label{RzetaTheta}
\pi^{-s/2}\Gamma(s/2)\Rzeta(s)=-\frac{e^{-\pi i
s/4}}{s}\int_{-1}^{-1+i\infty}\tau^{s/2}\vartheta_3'(\tau)\,d\tau.
\end{equation}
\end{theorem}
First, we study some simple properties of the theta functions involved.
We shall use the notation in \cite{WW}. Recall \cite{WW}*{p.~464} that $\vartheta_3(\tau)=\vartheta_3(0,\tau)$ is given by
\begin{equation*}
\vartheta_3(0,\tau)=1+2\sum_{n=1}^\infty e^{\pi i n^2\tau}.
\end{equation*}

We will be interested in the values of this function at the points
$\tau=-1+ix$
\begin{equation*}
\vartheta_3(0,-1+ix)=1+2\sum_{n=1}^\infty (-1)^n e^{-\pi n^2
x}=\vartheta_4(0, ix).
\end{equation*}
The transformation formula \cite{WW}*{p.~475--476},
\begin{equation*}
\vartheta_4(0,-\tau^{-1})=(-i\tau)^{1/2} \vartheta_2(0, \tau),
\end{equation*}
yields
\begin{equation*}
\vartheta_4(0, ix)=\frac{1}{\sqrt{x}}\vartheta_2(0, i/x)=
\frac{2}{\sqrt{x}}\sum_{n=0}^{\infty}e^{-\pi(n+\frac12)^2 /x}.
\end{equation*}

To alleviate the notation, we put
\begin{equation*}
\varphi(x):=1+2\sum_{n=1}^\infty (-1)^n e^{-\pi n^2 x}=
\frac{2}{\sqrt{x}}\sum_{n=0}^{\infty}e^{-\pi(n+\frac12)^2 /x},\qquad
x>0.
\end{equation*}

\begin{proposition}
We have the following inequalities
\begin{align*}
\varphi'(x)&\ge0,\qquad x>0,\\
\varphi'(x)&\le 2\pi e^{-\pi x},\qquad x>0,\\
\varphi'(x)&\le \frac{\pi}{2} x^{-5/2}e^{-\pi/4x},\qquad 0<x<1.
\end{align*}
\end{proposition}

\begin{proof}
From the product expansion of $\vartheta_4(0,ix)$ (see \cite{WW}*{p.~473})
\begin{displaymath}
\varphi(x)=\prod_{n=1}^\infty (1-e^{-\pi n x})(1-e^{-\pi(2n-1)x}),
\end{displaymath}
it follows that $\varphi(x)$ is increasing, therefore
$\varphi'(x)\ge0$. This is the first assertion in our Proposition.

We have
\begin{displaymath}
\varphi'(x)=2\pi e^{-\pi x} -2\pi\sum_{n=1}^\infty
\bigl\{(2n)^2e^{-\pi (2n)^2x}-(2n+1)^2e^{-\pi(2n+1)^2 x}\bigr\}.
\end{displaymath}
For $x>\frac{1}{4\pi}$ and $n\ge1$ we have
\begin{displaymath}
\bigl\{(2n)^2e^{-\pi (2n)^2x}-(2n+1)^2e^{-\pi(2n+1)^2
x}\bigr\}=\int_{2n}^{2n+1} 2t(\pi x t^2-1)e^{-\pi x t^2}\,dt\ge 0,
\end{displaymath}
since the integrand is positive. 
This proves the second assertion for $x>1/4\pi$. 

For the third assertion, observe that
\begin{displaymath}
\varphi'(x)
=\frac{\pi}{2}x^{-\frac{5}{2}}e^{-\frac{\pi}{4x}}-\sum_{n=1}^\infty
\Bigl\{x^{-\frac32} e^{-\frac{\pi}{x}(n-\frac12)^2}-2\pi\bigl(n+\tfrac12\bigr)^2x^{-\frac52}
e^{-\frac{\pi}{x}(n+\frac12)^2}\Bigr\}.
\end{displaymath}
For $0<x<1$, this is less than the first term because every term of
the sum is positive. To see it, we must show that
\begin{displaymath}
2\pi(n+\tfrac12)^2 x^{-1}e^{-2\pi n/x}\le 1, \qquad 0<x<1.
\end{displaymath}
The function $ye^{-y}$ is decreasing for $y>1$, thus $x^{-1}e^{-2\pi
n/x}\le e^{-2\pi n}$, and we only need to show that
\begin{displaymath}
2\pi(n+\tfrac12)^2 e^{-2\pi n}\le 1, \qquad 0<x<1.
\end{displaymath}
and this is easily checked for $n\ge1$.
This proves the third assertion. 

Now, if $0<x<1/4\pi$ we have
\begin{displaymath}
\varphi'(x)\le \frac{\pi}{2} x^{-5/2} e^{-\pi /4x}\le 1\le 2\pi e^{-\pi
x}.
\end{displaymath}
And the proposition is proved.
\end{proof}

\section{Infinite product.}

By means of  \eqref{RzetaTheta} we are now in a position to get bounds for $\Rzeta(s)$ on the hole plane.

\begin{proposition}\label{bound}
For $s=\sigma+it$
\begin{displaymath}
|s\pi^{-s/2}\Gamma(s/2) \Rzeta(s)|\le f(t) g(\sigma),
\end{displaymath}
where
\begin{displaymath}
f(t)=\begin{cases} 1& \text{if $t\ge0$,}\\
e^{\pi|t|/4}& \text{if $t<0$.}
\end{cases}
\qquad g(\sigma)=\begin{cases} 1& \text{if $\sigma\le 0$,}\\
3 \cdot2^{\sigma/4}\Gamma\left(1+\frac{\sigma}{2}\right)& \text{if
$\sigma> 0$.}
\end{cases}
\end{displaymath}
\end{proposition}

\begin{proof}
Taking $\tau=-1+ix$ in equation \eqref{RzetaTheta} we get
\begin{equation}\label{E:integralexpr}
s\,\pi^{-s/2}\Gamma(s/2)\Rzeta(s)=-e^{-\pi i s/4}\int_0^{+\infty}
(-1+ix)^{s/2} \varphi'(x)\,dx.
\end{equation}
(Observe that $\varphi'(x) =i\vartheta'_3(0,-1+ix)$).

It follows that
\begin{displaymath}
|s\pi^{-s/2}\Gamma(s/2)\Rzeta(s)|\le e^{\pi t/4}\int_0^{+\infty}
(1+x^2)^{\sigma/4}e^{-\frac{t}{2}(\pi-\arctan x)} \varphi'(x)\,dx.
\end{displaymath}
Now, if $t>0$
\begin{displaymath}
e^{\pi t/4}e^{-\frac{t}{2}(\pi-\arctan
x)}=e^{-\frac{t}{2}\left(\frac{\pi}{2}-\arctan x\right)}\le1
\end{displaymath}
and for $t\le 0$
\begin{displaymath}
e^{\pi t/4}e^{-\frac{t}{2}(\pi-\arctan x)}=e^{-\frac{\pi
t}{4}}e^{\frac{t}{2}\arctan x}\le e^{\pi |t|/4}.
\end{displaymath}

It follows that
\begin{displaymath}
|s\pi^{-s/2}\Gamma(s/2)\Rzeta(s)|\le f(t)\int_0^{+\infty}
(1+x^2)^{\sigma/4}\varphi'(x)\,dx.
\end{displaymath}

For $\sigma<0$ this is bounded by
$\int_0^{+\infty}\varphi'(x)\,dx\le \lim_{x\to+\infty}\varphi(x)=1$.
For $\sigma>0$ we have
\begin{multline*}
\int_0^{+\infty} (1+x^2)^{\sigma/4}\varphi'(x)\,dx\le
2^{\sigma/4}\int_0^1\varphi'(x)\,dx+\int_1^{+\infty}(2x^2)^{\sigma/4}\varphi'(x)\,dx\\
\le 2^{\sigma/4}+ 2\pi 2^{\sigma/4}\int_0^{+\infty}
x^{\sigma/2}e^{-\pi x}\,dx=
2^{\sigma/4}+2\pi^{-\sigma/2}2^{\sigma/4}\Gamma\left(1+\frac{\sigma}{2}\right)
\le 3\cdot 2^{\sigma/4}
\Gamma\left(1+\frac{\sigma}{2}\right).
\end{multline*}
The reader will not have difficulty proving the last inequality.
\end{proof}

\begin{theorem}
$\Rzeta(s)$ is an entire function of order $1$ and there exists a
constant $a$ such that
\begin{displaymath}
\Rzeta(s)=-\frac{e^{as}}{2\Gamma(1+\frac{s}{2})}\prod_{n\in\Z}
\Bigl(1-\frac{s}{\varrho_n}\Bigr)e^{s/\varrho_n}.
\end{displaymath}
Here $\Im\rho_n>0$ for $n\ge1$ and $\Im\rho_n\le 0$ for $n\le 0$. 
\end{theorem}

\begin{proof}
Consider the function $F(s)=s\pi^{-s/2}\Gamma(s/2)\Rzeta(s)$. By
\eqref{E:integralexpr} it is clear that it is an entire function. By
Proposition \ref{bound} we have for every $\varepsilon>0$ and $|s|>r_0$ 
\begin{equation}
|F(s)|\le 3e^{\pi |t|/4}\cdot
2^{|\sigma|/4}\Gamma\Bigl(1+\frac{|s|}{2}\Bigr)\le
C\exp\Bigl(\frac{1+\varepsilon}{2}|s|\log|s|\Bigr).
\end{equation}
Thus, $F(s)$ is an entire function of order $1$. Also, by Jensen's
Theorem, the number of zeros of $F(s)$ on the disc or center $0$ and
the radius $r$ are at most  $n(r)\le Cr\log r$.  From the Siegel results
\cite{Siegel} it follows that   $n(r)\ge c r\log r$.

By Theorem 7 in \cite{A98} we know that except for a region near the imaginary axis the only zeros of $\Rzeta(s)$ in the second and fourth quadrants are the trivial zeros at $-2$, $-4$, \dots\ These are not zeros of $F(s)$, since these zeros are simple \cite{A186} and the factor $\Gamma(1+s/2)$ have simple poles at these points. 

By Corollary 14 in \cite{A100}  there are  no zeros of $\Rzeta(s)$ for $s$ in the first quadrant for $\sigma\ge2$, except at most a finite number. 

Siegel proved that the number of zeros $\varrho=\beta+i\gamma$ of $\Rzeta(s)$ with $0<\gamma\le T$ is $\frac{T}{4\pi}\log\frac{T}{2\pi}-\frac{T}{4\pi}+o(T)$. We order, by increasing $\gamma$'s, these zeros of $F(s)$ in a sequence $\varrho_1$, $\varrho_2$, \dots, where each zero is repeated according to its multiplicity.

We have shown that there is a sequence of simple zeros of $\Rzeta(s)$ with $\gamma<0$ contained in the fourth quadrant \cite{A108}. In \cite{A100} Theorem 12 and Theorem 16 (and see also Remark 17 and \cite{A108}*{Prop.~1}) we see that all zeros of $\Rzeta(s)$ with $\gamma<0$ are contained in a narrow strip, which contains the zeros we have found in \cite{A108}. Possibly, these are the only ones. In any case, we order the zeros of $F(s)$ with $\gamma\le 0$ (there is a circle with center at $0$ in which, at present,  we cannot rule out the zeros of $\Rzeta(s)$ with $\gamma=0$ distinct from the trivial zeros) in an infinite sequence $\varrho_0$, $\varrho_{-1}$, $\varrho_{-2}$, \dots\

Since $n(r)\le Cr\log r$ we have $\sum_{n\in\Z}|\varrho_n|^{-2}<+\infty$.
Thus, we have
\begin{displaymath}
F(s)=e^{bs+c}\prod_{n\in\Z}\Bigl(1-\frac{s}{\varrho_n}\Bigr)e^{s/\varrho_n}.
\end{displaymath}

It follows that
\begin{displaymath}
\Rzeta(s)=\frac{\pi^{s/2}}{s\Gamma(s/2)}e^{bs+c}\prod_{n\in\Z}\Bigl(1-\frac{s}{\varrho_n}\Bigr)e^{s/\varrho_n}.
\end{displaymath}
It is known (see \cite{A166}) that $\Rzeta(0)=-\frac{1}{2}$. Hence, putting $s=0$ we
get
\begin{displaymath}
\Rzeta(s)=-\frac{e^{as}}{2\Gamma(1+s/2)}\prod_{n\in\Z}\Bigl(1-\frac{s}{\varrho_n}\Bigr)e^{s/\varrho_n}.
\end{displaymath}
\end{proof}

\begin{corollary}
\begin{equation}\label{E:rzetaDLog}
\frac{\Rzeta'(s)}{\Rzeta(s)}=a-\frac{1}{2}\frac{\Gamma'(1+s/2)}{\Gamma(1+s/2)}+\sum_{n\in\Z}
\Bigl(\frac{1}{s-\varrho_n}+\frac{1}{\varrho_n}\Bigr).
\end{equation}
\end{corollary}

We can  easily obtain the numerical value of the constant $a$,
putting  in the corollary $s=0$.
\begin{displaymath}
a=-2\Rzeta'(0)-\gamma/2.
\end{displaymath}
From which we obtain the numerical value
\begin{align*}
a=\; &{\scriptstyle
0.64087\,37327\,16376\,04132\,85108\,30622\,14339\,88326\,10314\,77196\,80921\,34404\,61888\,56880\,\dots}\\
{\scriptstyle+\; i\;}
&{\scriptstyle
0.55990\,02132\,94351\,56087\,92081\,78718\,19839\,03294\,83846\,48229\,05027\,43204\,55250\,68858\,\dots}
\end{align*}

\section{Phase of \texorpdfstring{$\Rzeta(s)$}{R(s)}.}
In \cite{AL}*{Proposition 2.4} it is proved that any real analytic function $f\colon\R\to\C$ can be written as $f(t)=U(t) e^{i\varphi(t)}$ where $U\colon\R\to\R$ and $\varphi\colon\R\to\R$ are real analytic. This decomposition is essentially unique, i.e.~if $U_1(t)e^{i\varphi_1(t)}=U_2(t)e^{i\varphi_2(t)}$, then we have either $U_1=U_2$ and $\varphi_1-\varphi_2=2k\pi$ or $U_1=-U_2$ and $\varphi_1-\varphi_2=(2k+1)\pi$. We say that $\varphi$ is the phase of $f$.

\begin{definition}
We will call $-\omega(t)$ the phase of $\Rzeta(\frac12+it)$ determined by 
\[\omega(0)=-\arg\Rzeta(\tfrac12)\approx2.86349.\]
Hence, $a\colon\R\to\R$ and $g\colon\R\to\R$ are real analytic functions such that
\[\Rzeta(\tfrac12+it)=e^{-i\;\omega(t)}g(t).\]
\end{definition}

\begin{remark}
It is possible that $\Rzeta(s)$ do not vanish on the critical line. In this case,
$g(t)=|\Rzeta(\frac12+it)|$  and $\omega(t)$ will be a value of $\arg\Rzeta(\frac12+it)$ for all $t\in\R$.
\end{remark}

The function $\omega(t)$ determine the zeros of $\zeta(s)$ on the critical line.

\begin{proposition}
The zeros of $\zeta(s)$ on the critical line are zeros of $\Rzeta(s)$ or 
the points where 
\[\cos(\vartheta(t)-\omega(t))=0.\]
\end{proposition}
\begin{proof}
In \cite{A166}*{eq.~(13) and (14)} it is shown that 
\[Z(t)=e^{i\vartheta(t)}\zeta(\tfrac12+it)=2\Re\bigl\{e^{i\vartheta(t)}\Rzeta(\tfrac12+it)\bigr\}=2\Re\bigl\{e^{i(\vartheta(t)-\omega(t))}g(t)\bigr\}.\]
From which the proposition follows.
\end{proof}

The zeros $\varrho=\beta+i\gamma$ of $\Rzeta(s)$ with $\gamma>0$ are connected to the zeros of $\zeta(s)$. Here, we define two counts of these zeros for $T>0$.
\begin{equation}\label{D:counts}
\begin{aligned}
N_r(T)&:=\#\bigl\{\varrho=\beta+i\gamma\colon\Rzeta(\varrho)=0,\quad  0<\gamma\le T, \quad\beta\ge\tfrac12\bigr\}\\
N_l(T)&:=\#\bigl\{\varrho=\beta+i\gamma\colon\Rzeta(\varrho)=0,\quad 0<\gamma\le T, \quad\beta\le\tfrac12\bigr\}
\end{aligned}
\end{equation}
Where $\#$ denotes the count of  zeros with its corresponding multiplicity, counting the zeros just on the critical line $\beta=1/2$ or at height $\gamma=T$ by one half. 
In this way $N_r(T)+N_l(T)$  is the number $N(T)$ of zeros of $\Rzeta(s)$ to height $T$.
Siegel \cite{Siegel} gives the main part of $N(T)$  and in \cite{A185} we have refined his Theorem showing
\begin{equation}\label{E:NT}
N(T)=\frac{T}{4\pi}\log\frac{T}{2\pi}-\frac{T}{4\pi}-\frac12 \sqrt{\frac{T}{2\pi}}+\Orden(T^{2/5}\log^2T)=\frac{\vartheta(T)-\pi\sqrt{T/2\pi}}{2\pi}+\Orden(T^{2/5}\log^2T).
\end{equation}

\begin{proposition}
For $T>0$ we have 
\begin{equation}\label{omegav}
\omega(T)=2\pi N_r(T)+\Orden(\log T).
\end{equation}
\end{proposition}

\begin{proof}
By \cite{A173}*{Proposition 1} for $\sigma\ge2$ and $t>32\pi$ we have $|\Rzeta(s)-1|<1$. 
Take $32\pi+1>t_0>32\pi$ such that $\Rzeta(s)$ do not have any zero on the horizontal line $t=t_0$. Take $T$ greater than  $t_0$ and such that there is no zero of $\Rzeta(s)$ with $\gamma=T$. We apply the argument principle to the rectangle $R=[1/2,2]\times[t_0,T]$.
\[N=\frac{1}{2\pi i}\int_\Gamma\frac{\Rzeta'(s)}{\Rzeta(s)}\,ds,\]
where $N$ is the number of zeros on the interior of $R$ and $\Gamma$ is the contour of $R$ with indentations on the possible zeros of $\Rzeta(s)$ on the border. 
By our elections, $\Rzeta(s)$ do not vanish on the contour of this rectangle, except at most a finite number of points on the border with $\sigma=\frac12$. We take a  semicircle with center at these possible zeros and radius $\varepsilon>0$ to exclude them from the interior of the  contour. 

The argument of $\Rzeta(s)$ changes only on a constant in the segment $t=t_0$. On the side $\sigma=2$ we have $|\Rzeta(s)-1|<1$ and the variation of the argument is at most $\pi$ on this segment.

On the horizontal side $\Im(s)=T$ the variation of the argument is $\le C\log T$. 
To see this, notice that by Proposition 6 in \cite{A173} for $T\ge T_0$ we have $|\Rzeta(2+iT)-1|\le 7/8$ so that $|\Rzeta(2+i T)|\ge 1/8$. In \cite{A92}*{Prop.~12} we prove that $|\Rzeta(\sigma+it)|\le2\sqrt{t/2\pi}$ for $\sigma>0$ and $t>16\pi$. Applying Backlund lemma \cite{A185}*{Lemma 1} to a disc with center at $2+iT$ and radius 2 we get
\[\Bigl|\Re\frac{1}{2\pi i}\int_{2+iT}^{1/2+iT}\frac{\Rzeta'(z)}{\Rzeta(z)}\,dz\Bigr|\le
\frac{1}{2}\Bigl(\log \frac{2\sqrt{(T+2)/2\pi}}{1/8}\Bigr)\Bigl(\log\frac{2}{3/2}\Bigr)^{-1}\ll \log T.\]

On the side $\sigma=1/2$ by \cite{AL}*{Prop. 2.10} we have
\[\omega(T)-\omega(t_0)=-\int_{t_0}^T\Re \frac{\Rzeta'(\frac12+ix)}{\Rzeta(\frac12+ix)}\,dx,\]
the integral being absolutely convergent. The integral on the corresponding side of $\Gamma$ when $\varepsilon\to0^+$ is equal to 
\[\lim_{\varepsilon\to0}\Re\frac{1}{2\pi i}\int_{L_\varepsilon}\frac{\Rzeta'(z)}{\Rzeta(z)}\,dz=\frac{\omega(T)-\omega(t_0)}{2\pi}-\frac{n}{2},\]
where $L_\varepsilon$ denotes the side of $\Gamma$ and $n$ the number of zeros of $\Rzeta(s)$ in the segment $[1/2+it_0,1/2+iT]$ counted with multiplicities. 

By the argument principle 
\[N=\frac{\omega(T)-\omega(t_0)}{2\pi}-\frac{n}{2}+\Orden(\log T).\]
Therefore $\omega(T)=2\pi(N+n/2)+\Orden(\log T)=2\pi N_r(T)+\Orden(\log T)$.
\end{proof}

Recall that the nontrivial zeros of $\Rzeta(s)$ repeated according to their multiplicity  are numerated by $\varrho_n=\beta_n+i\gamma_n$ 
for $n\in\Z$, in such a way that $\gamma_n>0$ for $n>0$ and $\gamma_n\le 0$ for $n\le 0$. We have also proved \cite{A100}*{Thm.~16} that for $n\le 0$ we have $\beta_n>0$ except at most for a finite number of zeros (our computation of zeros indicates that $\beta_n>0$ for all $n\le0$).

We need the following lemma
\begin{lemma}
We have
\[\sum_{n=1}^\infty\Bigl|\Re\frac{1}{\varrho_n}\Bigr|=\sum_{n=1}^\infty\frac{|\beta_n|}{\beta_n^2+\gamma_n^2}<+\infty, \qquad \sum_{n<0}\frac{\beta_n}{\beta_n^2+\gamma_n^2}=+\infty.\]
\end{lemma}

\begin{proof}
In the first case by Corollary 12 in  \cite{A98} and Corollary 14 in \cite{A100} we have for $n\ge n_0$ that $A\gamma_n^{2/5}\log\gamma_n>1-\beta_n$, and $\beta_n\le 2$. By \cite{A185}
\[n=\frac{\gamma_n}{4\pi}\log\frac{\gamma_n}{2\pi}-\frac{\gamma_n}{4\pi}+\Orden(\gamma^{1/2}_n).\]
Hence  $\frac{1}{5\pi}\gamma_n\log\gamma_n\le n\le \gamma_n\log\gamma_n$ for $n\ge n_0$, it follows that $\log\gamma_n\le \log n\le2\log\gamma_n$. Hence $\gamma_n\ge\frac{n}{\log\gamma_n}\ge\frac{n}{\log n}$. 

Then 
\begin{align*}
\sum_{n=1}^\infty\frac{|\beta_n|}{\beta_n^2+\gamma_n^2}&\le 
B+\sum_{n>n_0}\frac{A\gamma_n^{2/5}\log\gamma_n}{\gamma_n^2}\le B+A\sum_{n>n_0}\frac{\log n}{\gamma_n^{8/5}}\\ &\le B+A\sum_{n>n_0}\Bigl(\frac{\log n}{n}\Bigr)^{8/5}\log n<+\infty.\end{align*}

Except for a finite number the $\varrho_n$ with $n<0$ are those we have found in \cite{A108}. These zeros satisfy \cite{A108}*{Prop.~17}:
\[\beta_n\sim\frac{4\pi^2 |n|}{\log^2|n|},\quad \gamma_n\sim-\frac{4\pi |n|}{\log |n|}.\]
Therefore,
\[\sum_{n<0}\frac{|\beta_n|}{\beta_n^2+\gamma_n^2}\gg \sum_{n\le n'_0<0}
\frac{\beta_n}{2\gamma_n^2}\gg \sum_{n\le n'_0<0}
\frac{\frac{4\pi^2|n|}{\log^2|n|}}{2\frac{16\pi^2 n^2}{\log^2 |n|}}=\frac18\sum_{n\le n'_0<0}\frac{1}{|n|}=+\infty.\qedhere\]
\end{proof}

\begin{definition}
We define two functions $u\colon\R\to\R$ and $d\colon\R\to\R$ (for up and down)
by 
\begin{equation}\label{D:defud}
\begin{aligned}
u(t)&=\int_0^t\sum_{n>0}\frac{\beta_n-\frac12}{(\beta_n-\frac12)^2+(\gamma_n-x)^2}\,dx,\\
d(t)&=\int_0^t\sum_{n\le 0}\Bigl(\frac{\beta_n-\frac12}{(\beta_n-\frac12)^2+(\gamma_n-x)^2}-\frac{\beta_n-\frac12}{(\beta_n-\frac12)^2+\gamma_n^2}\Bigr)\,dx,
\end{aligned}
\end{equation}
\end{definition}

\begin{proposition}
There is a constant $B$ such that for $t\in\R$ we have
\begin{equation}\label{E:components}
\omega(t)=\vartheta(t)+u(t)+d(t)-Bt+\arctan(2t)+\omega(0).
\end{equation}
\end{proposition}
\begin{proof}
By  \cite{AL}*{Prop.~2.10} $\omega'(t)$  is 
\[\omega'(t)=-\Re\frac{\Rzeta'(\frac12+it)}{\Rzeta(\frac12+it)}.\]
Hence, by \eqref{E:rzetaDLog} $-\omega'(t)$ is the real part of 
\begin{displaymath}
\frac{\Rzeta'(\frac12+it)}{\Rzeta(\frac12+it)}=a-\frac{1}{\frac12+it}-\frac{1}{2}\frac{\Gamma'(\frac14+i\frac{t}{2})}{\Gamma(\frac14+i\frac{t}{2})}-\sum_n
\Bigl(\frac{1}{\beta_n-\frac12+i(\gamma_n-t)}-\frac{1}{\varrho_n}\Bigr).
\end{displaymath}
We prefer to put $\varrho_n-\frac12$ instead of $\varrho_n$, writing 
$\frac{\Rzeta'(\frac12+it)}{\Rzeta(\frac12+it)}$ as 
\begin{displaymath}
-\frac{1}{\frac12+it}-\frac{1}{2}\frac{\Gamma'(\frac14+i\frac{t}{2})}{\Gamma(\frac14+i\frac{t}{2})}-\sum_n
\Bigl(\frac{1}{\beta_n-\frac12+i(\gamma_n-t)}-\frac{1}{\beta_n-\frac12+i\gamma_n}\Bigr) + a+\sum_n\Bigl(\frac{1}{\varrho_n}-\frac{1}{\varrho_n-\frac12}\Bigr).
\end{displaymath}
Taking the real part
\begin{align*}
-\omega'(t)&=-\frac{\frac12}{\frac14+t^2}-\vartheta'(t)-\tfrac12\log\pi-\sum_{n\in\Z}
\Bigl(\frac{\beta_n-\frac12}{(\beta_n-\frac12)^2+(\gamma_n-t)^2}-
\frac{\beta_n-\frac12}{(\beta_n-\frac12)^2+\gamma_n^2}\Bigr)+C\\
&=-\frac{\frac12}{\frac14+t^2}-\vartheta'(t)-u'(t)-d'(t)+\sum_{n=1}^\infty\frac{\beta_n-\frac12}{(\beta_n-\frac12)^2+\gamma_n^2}+C-\frac12\log\pi.
\end{align*}
Hence $-\omega'(t)=-\frac{\frac12}{\frac14+t^2}-\vartheta'(t)-u'(t)-d'(t)+B$, where
\[B=-\frac12\log \pi+\sum_{n=1}^\infty\frac{\beta_n-\frac12}{(\beta_n-\frac12)^2+\gamma_n^2}+
\Re(a)+\Re\sum_n\Bigl(\frac{1}{\varrho_n}-\frac{1}{\varrho_n-\frac12}\Bigr).\]
Integrating the expression obtained for $\omega'(t)$ yields \eqref{E:components}.
\end{proof}

\begin{remark}
Using $t=0$ in the expression for $\frac{\Rzeta'(\frac12+it)}{\Rzeta(\frac12+it)}$ we get 
\[\frac{\Rzeta'(\frac12)}{\Rzeta(\frac12)}=-2-\frac{\Gamma'(1/4)}{2\Gamma(1/4)}+a+\sum_n\Bigl(\frac{1}{\varrho_n}-\frac{1}{\varrho_n-\frac12}\Bigr).\]
From which we obtain the numerical value
\[a+\sum_n\Bigl(\frac{1}{\varrho_n}-\frac{1}{\varrho_n-\frac12}\Bigr)=2+\frac{\Rzeta'(\frac12)}{\Rzeta(\frac12)}+\frac{\Gamma'(1/4)}{2\Gamma(1/4)}\]
\[=0.6373866805736784379 + i\;0.5524349167416397674.\]
We compute the partial sum 
\[\sum_{n=1}^{162215}\frac{\beta_n-\frac12}{(\beta_n-\frac12)^2+\gamma_n^2}=-0.009101619520742374473804923,\]
where the last summand is $-2.290228122019490122323182*10^{-11}$.
An approximate value for $B$ will be 
\[B\approx0.05592011\dots\]
\end{remark}

\begin{proposition}
The function $u(t)$ counts the difference between the number of zeros on the right and those on the left of the critical line. More precisely, for $t\to+\infty$ we have 
\begin{equation}\label{E:uvalue}
u(t)=\pi(N_r(t)-N_l(t))+\Orden(t^{1/2}\log t).
\end{equation}
\end{proposition}
\begin{proof}
By \eqref{D:defud} we have
\begin{equation}
u(t)=\int_0^t\sum_{n>0}\frac{\beta_n-\frac12}{(\beta_n-\frac12)^2+(\gamma_n-x)^2}\,dx=\sum_{\gamma>0,\beta\ne1/2}\Bigl(\arctan\frac{\gamma}{\beta-1/2}-\arctan\frac{\gamma-t}{\beta-1/2}\Bigr).
\end{equation}
The difference in arc tangents represents the angle subtended by the segment 
from $\frac12$ to $\frac12+it$ seen from the point $\varrho=\beta+i\gamma$, positive when the zero $\varrho$ is to the right of the critical line and negative when it is to the left. This is the sketch of our proof.

In this proof, it is convenient to denote by $\varrho=\beta+i\gamma$ a generic zero of $\Rzeta(s)$ with $\gamma>0$. Therefore, our $\gamma$ are always $>0$ and the sums refer to all the zeros of $\Rzeta(s)$ with $\gamma>0$ and $\beta\ne1/2$,  with the restrictions indicated in each case.

Let  $0<\theta<1$ be an exponent. In our reasoning, the error terms of order $t^\theta\log t$ and $t^{1-\theta}\log t$ will appear. Our best choice will be $\theta=1/2$.  But we will keep $\theta$ to make this clear. 

We divide the sum into four parts
\begin{equation}\label{E:foursums}
u(t)=\sum_{\gamma\le t^\theta}\cdots+\sum_{t^\theta<\gamma\le  t-t^\theta}\cdots+
\sum_{t-t^\theta<\gamma\le t+t^\theta}\cdots+
\sum_{t+t^\theta<\gamma}\Bigl(\arctan\frac{\gamma}{\beta-1/2}-\arctan\frac{\gamma-t}{\beta-1/2}\Bigr).\end{equation}
Let us call these sums $S_j$ for $1\le j\le 4$. We then estimate each $S_j$.
First, by \eqref{E:NT}
\begin{equation}\label{S_1}
|S_1|\le \pi N(t^\theta)\ll t^\theta\log t.
\end{equation}
Analogously,
\[|S_3|\le \pi (N(t+t^\theta)-N(t-t^\theta))=\pi\int_{t+t^\theta}^{t-t^\theta}
\Bigl(\frac{1}{4\pi}\log\frac{x}{2\pi}-\frac{1}{4\sqrt{2\pi x}}\Bigr)\,dx+\Orden(t^{2/5}\log^2 t).\]
It follows that (assuming $\theta>2/5$)
\begin{equation}\label{S_3}
|S_3|\ll t^\theta\log t+t^{2/5}\log^2 t\ll t^\theta\log t.
\end{equation}
$S_2$ is the main part of the sum. We divide this sum into two parts
\[S_2=S_2^r+S_2^l=\sum_{\substack{\beta>1/2\\t^\theta<\gamma\le  t-t^\theta}}\cdots
+\sum_{\substack{\beta<1/2\\t^\theta<\gamma\le  t-t^\theta}}\Bigl(\arctan\frac{\gamma}{\beta-1/2}-\arctan\frac{\gamma-t}{\beta-1/2}\Bigr),\]
notice that there are no terms with $\beta=1/2$.

We have (notice the change of $\gamma-t$ into $t-\gamma$)
\[S_2^r=\sum_{\substack{\beta>1/2\\t^\theta<\gamma\le  t-t^\theta}}
\Bigl(\arctan\frac{\gamma}{\beta-1/2}+\arctan\frac{t-\gamma}{\beta-1/2}\Bigr)\]
\[=\pi N_r(t^\theta,t-t^\theta)-
\sum_{\substack{\beta>1/2\\t^\theta<\gamma\le  t-t^\theta}}
\Bigl(\arctan\frac{\beta-1/2}{\gamma}+\arctan\frac{\beta-1/2}{t-\gamma}\Bigr).\]
Here $N_r(t^\theta,t-t^\theta)$ is the count of zeros with $t^\theta<\gamma\le t-t^\theta$ and $\beta>1/2$. For $t$ large enough we have $\beta-\frac12<2$ so that
\[\sum_{\substack{\beta>1/2\\t^\theta<\gamma\le  t-t^\theta}}
\Bigl(\arctan\frac{\beta-1/2}{\gamma}+\arctan\frac{\beta-1/2}{t-\gamma}\Bigr)
\ll \mskip-20mu\sum_{\substack{\beta>1/2\\t^\theta<\gamma\le  t-t^\theta}}
\Bigl(\frac{1}{\gamma}+\frac{1}{t-\gamma}\Bigr)\ll t^{-\theta}N(t)\ll t^{1-\theta}\log t.\]
\begin{equation}\label{S_2^r}
S_2^r=\pi N_r(t^\theta,t-t^\theta)+\Orden(t^{1-\theta}\log t).
\end{equation}
\begin{align*}
S_2^l&=\sum_{\substack{\beta<1/2\\t^\theta<\gamma\le  t-t^\theta}}\Bigl(\arctan\frac{\gamma}{\beta-1/2}-\arctan\frac{\gamma-t}{\beta-1/2}\Bigr)\\
&=-\sum_{\substack{\beta<1/2\\t^\theta<\gamma\le  t-t^\theta}}\Bigl(\arctan\frac{\gamma}{1/2-\beta}+\arctan\frac{t-\gamma}{1/2-\beta}\Bigr)\\
&=-\pi N_l(t^\theta,t-t^\theta)+\sum_{\substack{\beta<1/2\\t^\theta<\gamma\le  t-t^\theta}}\Bigl(\arctan\frac{1/2-\beta}{\gamma}+\arctan\frac{1/2-\beta}{t-\gamma}\Bigr).
\end{align*}
And we have 
\[\sum_{\substack{\beta<1/2\\t^\theta<\gamma\le  t-t^\theta}}\Bigl(\arctan\frac{1/2-\beta}{\gamma}+\arctan\frac{1/2-\beta}{t-\gamma}\Bigr)\ll
\sum_{\substack{\beta<1/2\\t^\theta<\gamma\le  t-t^\theta}}\Bigl(\frac{1/2-\beta}{\gamma}+\frac{1/2-\beta}{t-\gamma}\Bigr)\]
\[\ll t^{-\theta}\sum_{\substack{\beta<1/2\\t^\theta<\gamma\le  t-t^\theta}}(1/2-\beta)\ll t^{-\theta}\sum_{\substack{\beta<1/2\\t^\theta<\gamma\le  t-t^\theta}}(4-\beta)\ll t^{-\theta}\sum_{\beta<2, \gamma<t}(4-\beta).\]
By Siegel \cite{Siegel}*{eq.~(94)} (see also \cite{A102}*{eq.~(5)}) the last sum is $\ll t\log t$.
Therefore,
\begin{equation}\label{S_2^l}
S_2^l=-\pi N_l(t^\theta,t-t^\theta)+\Orden(t^{1-\theta}\log t).
\end{equation}
\[S_4=\sum_{t+t^\theta<\gamma}\Bigl(\arctan\frac{\gamma}{\beta-1/2}-\arctan\frac{\gamma-t}{\beta-1/2}\Bigr).\]
The two arctan have the same sign and then  $\arctan(a)-\arctan(b)=\arctan\frac{a-b}{1+ab}$ so that
\[|S_4|=\Bigl|\sum_{t+t^\theta<\gamma}\arctan\frac{t(\beta-1/2)}{(\beta-1/2)^2+\gamma(\gamma-t)}\Bigr|\le \sum_{t+t^\theta<\gamma}\arctan\frac{t|\beta-1/2|}{(\beta-1/2)^2+\gamma(\gamma-t)}.\]
Since $1-\beta\le A\gamma^{2/5}\log \gamma$ and $\gamma>t$
\[\frac{t|\beta-1/2|}{(\beta-1/2)^2+\gamma(\gamma-t)}\ll 
\frac{t\gamma^{2/5}\log \gamma}{\gamma(\gamma-t)}\ll t^{1-\theta}\frac{\gamma^{2/5}\log \gamma}{\gamma}=(t/\gamma)^{1-\theta}\gamma^{2/5-\theta}\log\gamma\le 1.\]
Hence,
\[|S_4|\ll\sum_{t+t^\theta<\gamma}\frac{t\gamma^{2/5}\log \gamma}{\gamma(\gamma-t)}
=\int_{t+t^\theta}^\infty\frac{tx^{2/5}\log x}{x(x-t)}\,dN(x)\]
\[\ll \int_{t+t^\theta}^\infty \frac{t x^{2/5}\log x}{x(x-t)}\Bigl(\frac{1}{4\pi}\log\frac{x}{2\pi}-\frac{1}{4\sqrt{2\pi x}}\Bigr)\,dx+
\int_{t+t^\theta}^\infty \frac{t x^{2/5}\log x}{x(x-t)}\,dR(x)\]
where $R(x)=\Orden (x^{2/5}\log^2 x)$. 
\begin{multline*}
\ll\int_{t+t^\theta}^\infty \frac{t x^{2/5}\log^2 x}{x(x-t)}\,dx+
\Bigl.\frac{t x^{2/5}\log x}{x(x-t)}R(x)\Bigr|_{x=t+t^\theta}^\infty\\-
\int_{t+t^\theta}^\infty R(x)\Bigl(\frac{t}{x^{8/5}(x-t)}-\frac{t\log x}{x^{3/5}(x-t)^2}-\frac{3t\log x}{5x^{8/5}(x-t)}\Bigr)\,dx.
\end{multline*}
We have 
\[\Bigl|\Bigl.\frac{t x^{2/5}\log x}{x(x-t)}R(x)\Bigr|_{x=t+t^\theta}^\infty\Bigr|=\Bigl|0-
\frac{t (t+t^\theta)^{2/5}\log(t+t^\theta)}{(t+t^\theta)t^\theta}R(t+t^\theta)\Bigr|\ll t^{4/5-\theta}\log^3t.\]
Therefore,
\begin{multline*}
|S_4|\ll \Orden(t^{4/5-\theta}\log^3t)+\int_{t+t^\theta}^\infty \frac{t \log^2 x}{x^{3/5}(x-t)}\,dx
+\int_{t+t^\theta}^\infty \frac{t\log^2 x}{x^{6/5}(x-t)}\,dx\\+
\int_{t+t^\theta}^\infty\frac{t\log^3 x}{x^{1/5}(x-t)^2}\,dx+
\int_{t+t^\theta}^\infty\frac{t\log^3 x}{x^{6/5}(x-t)}\,dx.
\end{multline*}
Some of these integrals dominate the others, so that we obtain 
\[|S_4|\ll \Orden(t^{4/5-\theta}\log^3t)+\int_{t+t^\theta}^\infty \frac{t \log^2 x}{x^{3/5}(x-t)}\,dx
+\int_{t+t^\theta}^\infty\frac{t\log^3 x}{x^{1/5}(x-t)^2}\,dx.
\]
From which it is easy to get 
\begin{A}
\[\int_{t+t^\theta}^\infty \frac{t \log^2 x}{x^{3/5}(x-t)}\,dx=t^{2/5}\int_{1+t^{\theta-1}}^\infty\frac{(\log y+\log t)^2\,dy}{y^{3/5}(y-1)}\]
\[\ll t^{2/5}\log^2t+t^{2/5}\int_{1+t^{\theta-1}}^2\frac{(\log y+\log t)^2\,dy}{y^{3/5}(y-1)}\]
\[\ll t^{2/5}\log^2t+t^{2/5}\log^2t\int_{1+t^{\theta-1}}^2\frac{dy}{y-1}+
t^{2/5}\log t\ll t^{2/5}\log^3t.\]
\[\int_{t+t^\theta}^\infty\frac{t\log^3 x}{x^{1/5}(x-t)^2}\,dx=
t^{-1/5}\int_{1+t^{\theta-1}}^\infty\frac{(\log y+\log t)^3}{y^{1/5}(y-1)^2}\,dy\]
\[\ll t^{-1/5}\log^3t+t^{-1/5}\int_{1+t^{\theta-1}}^2\frac{(\log y+\log t)^3}{y^{1/5}(y-1)^2}\,dy\ll t^{-1/5}\log^3t+t^{-1/5}\log^3t\int_{1+t^{\theta-1}}^2\frac{dy}{(y-1)^2}\]
\[\ll t^{4/5-\theta}\log^3t.\]
\end{A}
\begin{equation}\label{S_4}
|S_4|=\Orden(t^{4/5-\theta}\log^3t)+\Orden(t^{2/5}\log^3t)+\Orden(t^{1-\theta-1/5}\log^3t)=\Orden(t^{2/5}\log^3t).
\end{equation}
The main part of $S_2^r+S_2^l$ is $ \pi (N_r(t^\theta,t-t^\theta)-N_l(t^\theta,t-t^\theta)$. This is equal to 
\[\pi (N_r(t-t^\theta)-N_l(t-t^\theta))-\pi (N_r(t^\theta)-N_l(t^\theta)),\]
since the zeros with $\beta=1/2$ are counted in the same form in $N_t(t)$ and $N_l(t)$ so that in the difference they are not counted. 
We have 
\[|(N_r(t-t^\theta)-N_l(t-t^\theta))- (N_r(t^\theta)-N_l(t^\theta))-(N_r(t)-N_l(t))|\]
\[\le |N_r(t-t^\theta)-N_r(t)|+|N_l(t-t^\theta)-N_l(t)|+|N_r(t^\theta)|+|N_l(t^\theta)|\le 2|N(t-t^\theta)-N(t)|+2N(t^\theta)\]
\[\ll t^\theta\log t+t^\theta\log t\ll t^\theta\log t.\]
From this and from \eqref{S_2^r} and \eqref{S_2^l} it follows that 
\begin{equation}\label{S_2}
S_2=\pi(N_r(t)-N_l(t))+\Orden(t^\theta\log t)+\Orden(t^{1-\theta}\log t).
\end{equation}
Combining \eqref{E:foursums}, \eqref{S_1}, \eqref{S_2}, \eqref{S_3} and \eqref{S_4}, and taking $\theta=1/2$ we obtain \eqref{E:uvalue}.
\end{proof}

\begin{proposition}
For $t>0$ we have 
\begin{equation}\label{dvalue}
d(t)=Bt-\tfrac12\vartheta(t)+\Orden(t^{1/2}\log t).
\end{equation}
Then we derive 
\begin{equation}
\omega(t)=\tfrac12\vartheta(t)+u(t)+\Orden(t^{1/2}\log t).
\end{equation}
\end{proposition}
\begin{proof}
In equation \eqref{E:components}
\[\omega(t)=\vartheta(t)+u(t)+d(t)-Bt+\arctan(2t)+\omega(0),\]
substitute $\omega(t)$ by its value \eqref{omegav}, $u(t)$ by \eqref{E:uvalue}. This yields
\[2\pi N_r(t)+\Orden(\log t)=\vartheta(t)+\pi(N_r(t)-N_l(t))+\Orden(t^{1/2}\log t)+d(t)-Bt+\arctan(2t)+\omega(0)\]
Or equivalently
\[\pi(N_r(t)+N_l(t))=\vartheta(t)+d(t)-Bt+\Orden(t^{1/2}\log t).\]
From \eqref{E:NT} we get  $N_r(t)+N_l(t)=\frac{1}{2\pi}\vartheta(t)+\Orden(t^{1/2})$. Substituting this value into the above equation, we get \eqref{dvalue}. The last equation for $\omega(t)$ follows by substituting the value of $d(t)$ into \eqref{E:components}. 
\end{proof}

\begin{remark}
Our zero computations \cite{A172} are consistent with $N_r(t)\approx\frac12N_l(t)$. If this were true, we would have
\begin{equation}
N_l(t)\approx \frac{\vartheta(t)}{3\pi},\quad
N_r(t)\approx \frac{\vartheta(t)}{6\pi},\quad
u(t)\approx -\frac{\vartheta(t)}{6},\quad
d(t)\approx -\frac{\vartheta(t)}{2}+Bt,\quad 
\omega(t)\approx\frac{\vartheta(t)}{3}.
\end{equation}
\end{remark}


\begin{thebibliography}{999}

\bibitem{AL}
\textsc{J. Arias de Reyna \&\ J. van de Lune}, \emph{On the exact location of the non-trivial zeros of Riemann's zeta function}, Acta Arith. \textbf{163} (2014) 215--245. \href{https://www.impan.pl/en/publishing-house/journals-and-series/acta-arithmetica/all/163/3/83430/on-the-exact-location-of-the-non-trivial-zeros-of-riemann-s-zeta-function}{DOI: 10.4064/aa163-3-3}.

\bibitem{A166}
\textsc{Arias de Reyna, J.},  \emph{Riemann's auxiliary function. Basic results}, \href{https://arxiv.org/abs/2406.02403}{arXiv:2406.02403}. 

\bibitem{A172}
\textsc{Arias de Reyna, J.},  \emph{Statistics of zeros of the auxiliary function}, \href{https://arxiv.org/abs/2406.03041}{arXiv:2406.03041}. 



\bibitem{A92}
\textsc{Arias de Reyna, J.},  \emph{Simple bounds of zeta and related functions}, preprint (92)

\bibitem{A98}
\textsc{Arias de Reyna, J.},  \emph{Regions without zeros for the auxiliary function of Riemann}, \href{https://arxiv.org/abs/2406.03825}{arXiv:} \href{https://arxiv.org/abs/2406.03825}{2406.03825}. 




\bibitem{A100}
\textsc{Arias de Reyna, J.},  \emph{Asymptotic expansions of the auxiliary function}, \href{https://arxiv.org/abs/2406.04714}{arXiv:2406.04714}. 


\bibitem{A173}
\textsc{Arias de Reyna, J.},  \emph{Riemann's auxiliary function. Right limit of zeros}, \href{https://arxiv.org/abs/2406.07014}{arXiv:2406.07014}. 

\bibitem{A102}
\textsc{Arias de Reyna, J.},  \emph{On Siegel results about the zeros of the auxiliary function of Riemann}, \href{https://arxiv.org/abs/2406.07968}{arXiv:2406.07968}. 


\bibitem{A185}
\textsc{Arias de Reyna, J.},  \emph{On the number of zeros of $\Rzeta(s)$}, \href{https://arxiv.org/abs/2406.08890}{arXiv:2406.08890}. 

\bibitem{A186}
\textsc{Arias de Reyna, J.},  \emph{Trivial zeros of Riemann auxiliary function}, \href{https://arxiv.org/abs/2406.09796}{arXiv:2406.09796}. 


\bibitem{A108}
\textsc{Arias de Reyna, J.},  \emph{Zeros of $\Rzeta(s)$ in the fourth quadrant},  \href{https://arxiv.org/abs/2406.11279}{arXiv:2406.11279}. 




\bibitem{Con}
\textsc{Conrey J.B.}, \emph{More than two fifths of the zeros of the Riemann zeta function are on the critical line}, J. Reine Angew. Math. \textbf{399} (1989) 1--26.



\bibitem{Lev}
\textsc{Levinson N.}, \emph{More than one third of zeros of Riemann’s zeta-function are on $\sigma=\frac12$}, Adv. Math. \textbf{13} (1974) 383--436.




\bibitem{F}
\textsc{Shaoji Feng}, \emph{Zeros of the Riemann zeta function on the critical line}, J. Number Theory, \textbf{132}, (2012), 511--542.

\bibitem{P}
\textsc{Pratt, K.;  Robles, N.; Zaharescu, A.; Zeindler, D.},
\emph{More than five-twelfths of the zeros of $\zeta$ are on the critical line},
Res. Math. Sci. \textbf{7} (2020), no. 2, Paper No. 2, 74 pp.

\bibitem{Siegel}
\textsc{C. L. Siegel}, \emph{Uber Riemann's Nachla{\ss} zur
analytischen Zahlentheorie}, {\rm Quellen und Studien zur Geschichte
der Mathematik, Astronomie und Physik} \textbf{2} (1932), 45--80.
Reprinted in~\cite{SW}, 1, 275--310.  \href{https://arxiv.org/abs/1810.05198}{English version}.

\bibitem{SW}
\textsc{C. L. Siegel}, \emph{Carl Ludwig Siegel's Gesammelte Abhandlungen}, 
(edited by K. Chandrasekharan and H. Maa\ss), Springer-Verlag, Berlin, 1966.


\bibitem{WW} 
\textsc{E.\ T.\ Whittaker, G. N. Watson}, \emph{A Course of Modern Analysis}, Cambridge University Press,  1965.


\end{thebibliography}
\end{document}